\documentclass[aos,MSNbibl,nameyear,dvips]{arximspdf}
\usepackage{graphicx}

\doi{10.1214/10-AOS870}
\volume{39}
\issue{2}
\pubyear{2011}
\firstpage{1241}
\lastpage{1265}

\makeatletter

\def\bptnote#1{}

\def\bfx{\mathbf{x}}
\def\bfX{\mathbf{X}}
\def\bfl{\mathbf{l}}
\def\bfA{\mathbf{A}}
\def\bfC{\mathbf{C}}
\def\bfS{\mathbf{S}}
\def\bfI{\mathbf{I}}
\def\bfmu{{\bolds{\mu}}}
\def\varep{\varepsilon}
\def\bfdelta{{\bolds{\delta}}}
\def\bfxi{\bolds{\xi}}
\def\bfSigma{{\bolds{\Sigma}}}
\def\RB{R_{\mathrm{OPT}}}
\def\bfzero{\bolds{0}}
\def\zetap{\Delta_p^2}
\def\zetaps{\Delta_p }
\def\rightarrowp{\rightarrow_{P}}

\newproclaim{remark}{Remark}
\newproclaim{definition}{Definition}
\newtheorem{lemma}{Lemma}

\makeatother

\begin{document}
\begin{frontmatter}

\title{Sparse linear discriminant analysis by thresholding for high dimensional data}
\runtitle{Sparse linear discriminant analysis}

\begin{aug}
\author[A]{\fnms{Jun} \snm{Shao}\corref{}\ead[label=e1]{shao@stat.wisc.edu}\thanksref{aut1}},
\author[A]{\fnms{Yazhen} \snm{Wang}\ead[label=e2]{yzwang@stat.wisc.edu}\thanksref{aut2}},
\author[A]{\fnms{Xinwei} \snm{Deng}\ead[label=e3]{xdeng@stat.wisc.edu}}
and
\author[A]{\fnms{Sijian} \snm{Wang}\ead[label=e4]{wangs@stat.wisc.edu}}
\thankstext{aut1}{Supported in part by the NSF Grant SES-0705033.}
\thankstext{aut2}{Supported in part by the NSF Grant DMS-10-05635.}
\runauthor{Shao, Wang, Deng and Wang}
\affiliation{East China Normal University and University of Wisconsin}
\address[A]{Department of Statistics\\
University of Wisconsin\\
1300 University Ave.\\
Madison, Wisconsin 53706\\
USA \\
\printead{e1}\\
\phantom{E-mail:\ }\printead*{e2}\\
\phantom{E-mail:\ }\printead*{e3}\\
\phantom{E-mail:\ }\printead*{e4}} 
\end{aug}

\received{\smonth{2} \syear{2010}}
\revised{\smonth{9} \syear{2010}}

%
\begin{abstract}
In many social, economical, biological and medical studies,
one objective is to classify a subject into one of
several classes based on a set of variables observed from
the subject.
Because the probability distribution of the variables is
usually unknown, the rule of classification is constructed
using a training sample.
The well-known linear discriminant analysis (LDA) works well
for the situation where the number of variables used for
classification is much smaller than the training sample size.
Because of the advance in technologies,
modern statistical studies often face classification problems
with the number of variables much larger than
the sample size, and the LDA may perform poorly.
We explore when and why the LDA has poor performance and
propose a sparse LDA that is
asymptotically optimal
under some sparsity conditions on the unknown parameters.
For illustration of application, we
discuss an example of classifying human cancer into two classes
of leukemia based on a set of 7,129 genes and
a training sample of size 72.
A simulation is also conducted to check the
performance of the proposed method.
\end{abstract}

%


\end{frontmatter}

\section{Introduction}\label{sec1}

The objective of a classification problem is to classify
a subject to one of several classes based on a
$p$-dimensional vector $\mathbf{x}$ of characteristics
observed from the subject. In most applications, variability exists,
and hence
$\mathbf{x}$ is random. If the distribution of $\mathbf{x}$ is
known, then we can construct an optimal classification rule that
has the smallest possible misclassification rate.
However, the distribution of $\mathbf{x}$ is usually unknown, and
a~classification rule has to be constructed using a training sample.
A~statistical issue is how to use the training sample to
construct a classification rule that has a misclassification rate
close to that of the optimal rule.

In traditional applications, the dimension $p$ of $\mathbf{x}$ is fixed
while the training sample size $n$ is large.
Because of the advance in technologies,
nowadays a~much larger
amount of information can be collected, and the resulting
$\mathbf{x}$ is of a high dimension. In many recent applications, $p$
is much larger than the training sample size, which
is referred to as the large-$p$-small-$n$ problem
or ultra-high dimension problem when $p= O(e^{n^\beta})$
for some $\beta\in(0,1)$.
An example is a study with genetic or microarray data.
In our example presented in Section \ref{sec5}, for instance, a
crucial step for a successful chemotherapy treatment is to
classify human cancer into two classes of leukemia, acute myeloid
leukemia and
acute lymphoblastic leukemia, based on $p=7\mbox{,}129$
genes and a~training sample
of 72 patients. Other examples include data from radio\-logy, biomedical imaging,
signal processing, climate and finance.
Although more information is better when the distribution of
$\mathbf{x}$ is known, a larger dimension $p$ produces more
uncertainty when the distribution of $\mathbf{x}$ is unknown and,
hence, results in a greater challenge for data analysis
since the training sample size $n$ cannot increase as fast as $p$.

The well-known linear discriminant analysis (LDA) works well
for fixed-$p$-large-$n$ situations and is asymptotically
optimal in the sense that, when $n$ increases to infinity,
its misclassification rate over that of the optimal rule
converges to one. In fact, we show in this
paper that the LDA is still asymptotically optimal
when $p$ diverges to
infinity at a rate slower than $\sqrt{n}$. On the other hand,
\citet{BicLev04} showed that the LDA is asymptotically
as bad as random guessing when $p>n$; some similar results are
also given in this paper. The main purpose of this paper
is to construct a sparse LDA and show it is asymptotically
optimal under some sparsity conditions
on unknown parameters and some condition on the
divergence rate of $p$
(e.g., $n^{-1} \log p \rightarrow0$ as $n \rightarrow\infty$).
Our proposed sparse LDA is based on the thresholding
methodology, which was developed in wavelet
shrinkage for function estimation [\citet{DonJoh94},
Donoho et al. (\citeyear{Donetal95})] and covariance matrix
estimation [\citet{BicLev08}].
There exist a few other sparse LDA methods, for example,
\citet{GuoHasTib07}, \citet{CleHas} and
\citet{QiaZhoHua09}. The key differences between the existing methods
and ours are the conditions on sparsity and
the construction of sparse estimators of parameters.
However, no asymptotic results
were established in the existing papers.

For high-dimensional $\mathbf{x}$ in regression,
there exist some variable selection
methods [see a recent review by \citet{FanLv10}].
For constructing a~classification rule using variable selection,
we must identify not only
components of $\mathbf{x}$ having mean effects for classification,
but also components of $\mathbf{x}$ having effects
for classification through their correlations with
other components [see, e.g., \citet{KohJoh97},
\citet{ZhaWan10}].
This may be a~very difficult task when $p$ is much
larger than $n$, such as $p=7\mbox{,}129$ and $n=72$ in the
leukemia example in Section \ref{sec5}. Ignoring the correlation,
\citet{FanFan08} proposed the features annealed
independence rule (FAIR), which first selects
$m$ components of $\mathbf{x}$ having mean effects for
classification and then
applies the naive Bayes rule (obtained by assuming that components of
$\mathbf{x}$ are independent)
using the selected $m$ components of $\mathbf{x}$ only.
Although no sparsity condition on the covariance matrix of $\mathbf
{x}$ is
required,
the FAIR is not asymptotically optimal
because the correlation between 
components of $\mathbf{x}$ is ignored.
Our approach is not a variable selection approach,
that is, we do not try to identify a subset of components
of $\mathbf{x}$ with a size smaller than $n$.
We use thresholding estimators of the mean effects
as well as Bickel and Levina's (\citeyear{BicLev08})
thresholding estimator of the covariance matrix
of $\mathbf{x}$, but we allow the number of nonzero estimators (for
the mean differences or covariances) to be much larger
than $n$ to ensure the asymptotic optimality of the
resulting classification rule.

The rest of this paper is organized as follows.
In Section \ref{sec2}, after introducing some notation and
terminology, we establish a sufficient condition on the
divergence of $p$ under which the LDA is still asymptotically
close to the optimal rule. We also show that,
when $p$ is large compared with $n$ ($p/n \rightarrow
\infty$), the performance of the LDA is not good
even if we know the covariance matrix of $\mathbf{x}$,
which indicates the need of sparse estimators
for both the mean difference and covariance matrix.
Our main result is given in Section \ref{sec3}, along with
some discussions about various sparsity conditions
and divergence rates of $p$ for which the proposed
sparse LDA performs well asymptotically.
Extensions of the main result are discussed in Section \ref{sec4}.
In Section \ref{sec5}, the proposed sparse LDA
is illustrated in the example
of classifying human cancer into two classes of leukemia, along with
some simulation results for examining misclassification rates.
All technical proofs are given in Section \ref{sec6}.

\section{The optimal rule and linear discriminant analysis}\label{sec2}

We focus on the classification problem with two classes. The general
case with three or more classes is discussed in Section \ref{sec4}.
Let $\mathbf{x}$ be a $p$-dimensional normal random vector belonging to
class $k$ if $\mathbf{x}\sim N_p( \bfmu_k , \bfSigma)$, $k=1,2$,
where $ \bfmu_1 \neq\bfmu_2$, and $\bfSigma$ is positive definite.
The misclassification rate of any classification rule is the average
of the probabilities of making two types of misclassification:
classifying $\mathbf{x}$ to class 1 when $\mathbf{x}\sim N_p(\bfmu_2 ,
\bfSigma)$
and classifying $\mathbf{x}$ to class 2 when $\mathbf{x}\sim
N_p(\bfmu_1 ,
\bfSigma)$.

If $\bfmu_1$, $\bfmu_2$ and $\bfSigma$ are known, then
the optimal classification rule, that is, the rule with the smallest
misclassification
rate, classifies $\mathbf{x}$ to class 1 if and only if
$\bfdelta' \bfSigma^{-1} ( \mathbf{x}- \bar{\bfmu}) \geq0$,
where $\bar\bfmu= (\bfmu_1 + \bfmu_2)/2$,
$\bfdelta= \bfmu_1- \bfmu_2$, and
$\mathbf{a}'$ denotes the transpose of the vector $\mathbf{a}$.
This rule is also the Bayes rule with equal prior probabilities for two classes.
Let $\RB$ denote the misclassification rate of the optimal rule.
Using the normal distribution, we can show that
%
\begin{equation}\label{delta}
\RB= \Phi( - \zetaps/2 ),\qquad
\zetaps= \sqrt{\bfdelta' \bfSigma^{-1} \bfdelta},
\end{equation}
where $\Phi$ is the standard normal distribution function.
Although $0 < \RB< 1/2$, $\RB\rightarrow0$ if $\zetaps
\rightarrow\infty$ as $p \rightarrow\infty$
and $\RB\rightarrow1/2$ if $\zetaps\rightarrow0$.
Since $1/2$ is the misclassification rate of random guessing,
we assume the following regularity conditions:
there is a constant $c_0$ (not depending on $p$) such that
%
\begin{equation} \label{conds}
\mbox{$c_0^{-1} \leq$ all eigenvalues of $\bfSigma\leq c_0$}
\end{equation}
and
%
\begin{equation}\label{condd}
c_0^{-1} \leq\max_{j \leq p} \delta_j^2 \leq c_0 ,
\end{equation}
where $\delta_j$ is the $j$th component of $\bfdelta$.
Under (\ref{conds})--(\ref{condd}), $\zetaps\geq c_0^{-1}$, and hence
$\RB\leq\Phi( - (2c_0)^{-1})< 1/2$.
Also, $\zetap= O( \| \bfdelta\|^2 )$ and
$\| \bfdelta\|^2 = O(\zetap)$ so that
the rate of $\| \bfdelta\|^2 \rightarrow\infty$
is the same as the rate of $\zetap\rightarrow\infty$,
where $\| \mathbf{a}\|$ is the $L_2$-norm of the vector $\mathbf{a}$.

In practice, $\bfmu_k$ and $\bfSigma$ are typically unknown, and we
have a
training sample $\bfX= \{ \mathbf{x}_{ki}, i=1,\ldots,n_k, k=1,2 \}
$, where
$n_k$ is the sample size for class~$k$,
$\mathbf{x}_{ki} \sim N_p( \bfmu_k , \bfSigma)$, $k=1,2$,
all $\mathbf{x}_{ki}$'s are independent and
$\bfX$ is independent of $\mathbf{x}$ to be classified.
The limiting process considered in this paper is
the one with $n = n_1+n_2 \rightarrow\infty$.
We assume that $n_1/n $ converges to a constant strictly between 0 and 1;
$p$ is a function of $n$, but the subscript $n$ is omitted for simplicity.
When $n \rightarrow\infty$, $p$ may diverge to $\infty$,
and the limit of $p/n$ may be 0, a positive constant, or $\infty$.

For a classification rule $T$ constructed using the training sample,
its performance can be assessed by the conditional misclassification
rate $R_T (\bfX)$
defined as the average of the conditional probabilities of
making two types of misclassification,
where the conditional probabilities are with respect to~$\mathbf{x}$,
given the training sample~$\bfX$.
The unconditional misclassification rate is
$R_T = E[ R_T(\bfX)]$. The asymptotic performance of $T$
refers to the limiting behavior of $R_T(\bfX) $ or $R_T$ as $n
\rightarrow\infty$.
Since $0 \leq R_T( \bfX) \leq1$, by the dominated convergence theorem,
if $R_T( \bfX)\rightarrowp c$, where $c$ is a constant and
$\rightarrowp$ denotes
convergence in probability, then $R_T \rightarrow c$.
Hence, in this paper we focus on the limiting behavior of the
conditional misclassification
rate $R_T( \bfX)$.

We hope to find a rule $T$ such that $R_T (\bfX)$ converges in probability
to the same limit as $\RB$, the misclassification rate of the
optimal rule.
If $\RB\rightarrow0$, however, we hope
not only $R_T (\bfX) \rightarrowp0$, but also
$R_T(\bfX)$ and $\RB$ have the same convergence rate.
This leads to the following definition.

\begin{definition}\label{defin1}
Let $T$ be a classification rule with conditional
misclassification rate $R_T(\bfX)$, given the training sample $\bfX$.
\begin{longlist}[(iii)]
\item[(i)] $T$ is asymptotically optimal if
$ R_T (\bfX)/ \RB\rightarrowp1 $.
\item[(ii)] $T$ is asymptotically sub-optimal if
$R_T(\bfX) - \RB\rightarrowp0 $.
\item[(iii)] $T$ is asymptotically worst if
$R_T (\bfX) \rightarrowp1/2 $.
\end{longlist}
\end{definition}

If $\lim_{n \rightarrow\infty} \RB>0$ [i.e., $\zetaps$ in (\ref
{delta}) is
bounded], then the asymptotic sub-optimality is the same as
the asymptotic optimality.
Part (iii) of Definition \ref{defin1} comes from the fact that $1/2$ is the
misclassification rate
of random guessing.

In this paper we focus on the classification rules of the form
%
\begin{equation}
\mbox{classifying $\mathbf{x}$ to class 1 if and only if }
\hat\bfdelta{}' \hat\bfSigma{}^{-1} ( \mathbf{x}- \hat{\bar\bfmu} )
\geq0, \label{rule}
\end{equation}
where $\hat\bfdelta$, $\hat{\bar\bfmu}$ and $\hat\bfSigma{}^{-1}$
are estimators of $\bfdelta$, $\bar\bfmu$ and $\bfSigma^{-1}$, respectively,
constructed using the training sample $\bfX$.

The well-known linear discriminant analysis (LDA) uses
the maximum likelihood estimators $\bar{\mathbf{x}}_1$, $\bar
{\mathbf{x}}_2$
and $\bfS$, where
\[
\bar{\mathbf{x}}_k = \frac{1}{n_k} \sum_{i=1}^{n_k} \mathbf
{x}_{ki} , \qquad
k=1,2,\qquad
\bfS= \frac{1}{n} \sum_{k=1}^2
\sum_{i=1}^{n_k} ( \mathbf{x}_{ki}- \bar{\mathbf{x}}_k)( \mathbf
{x}_{ki}- \bar
\mathbf{x}_k)' .
\]
The LDA is given by (\ref{rule}) with
$\hat\bfdelta= \bar{\mathbf{x}}_1 - \bar{\mathbf{x}}_2$, $\hat
{\bar\bfmu} =
\bar{\mathbf{x}}= (\bar{\mathbf{x}}_1+\bar{\mathbf{x}}_2)/2$,
$\hat\bfSigma{}^{-1} =
\bfS^{-1}$ when $\bfS^{-1}$ exists, and $\hat\bfSigma{}^{-1} = $
a generalized inverse $\bfS^-$ when $\bfS^{-1}$ does not exist (e.g.,
when $p> n$).
A straightforward calculation shows that, given
$\bfX$, the conditional misclassification rate of the LDA is
%
\begin{equation}\label{rate}
\frac{1}{2} \sum_{k=1}^2 \Phi\biggl(
\frac{(-1)^k \hat\bfdelta{}' \hat\bfSigma{}^{-1} ( \bfmu_k -\bar
\mathbf{x}_k )
- \hat\bfdelta{}' \hat\bfSigma{}^{-1} \hat\bfdelta/2}{\sqrt{\hat
\bfdelta{}'\bfS^{-1} \bfSigma\hat\bfSigma{}^{-1} \hat\bfdelta}} \biggr).
\end{equation}

Is the LDA asymptotically optimal or sub-optimal according to
Definition~\ref{defin1}?
Bickel and Levina [(\citeyear{BicLev04}), Theorem 1]
showed that, if $p > n$ and $p / n \rightarrow\infty$, then
the unconditional misclassification rate of the LDA converges to $1/2$
so that the LDA is asymptotically worst.
A natural question is, for what kind of $p $ (which may diverge to
$\infty$), is
the LDA asymptotically optimal or sub-optimal.
The following result provides an answer.

\begin{thm}\label{thm1}
Suppose that (\ref{conds})--(\ref{condd}) hold
and $s_n = p \sqrt{\log p}/\sqrt{n} \rightarrow0$.
\begin{longlist}[(iii)]
\item[(i)] The conditional misclassification rate of the LDA
is equal to
\[
R_{\mathrm{LDA}} (\bfX) = \Phi\bigl( - [1+O_P(s_n)]\zetaps/2 \bigr).
\]
\item[(ii)] If $\zetaps$ is bounded, then the LDA is asymptotically
optimal and
\[
\frac{R_{\mathrm{LDA}} ( \bfX) }{\RB} - 1 = O_P(s_n) .
\]
\item[(iii)] If $\zetaps\rightarrow\infty$, then the LDA is
asymptotically sub-optimal.

\item[(iv)] If $\zetaps\rightarrow\infty$ and
$s_n \zetap= ( p \sqrt{\log p }/\sqrt{n}) \zetap\rightarrow0 $,
then the LDA is asymptotically optimal.
\end{longlist}
\end{thm}

\begin{remark}\label{rem1}
Since $\zetaps\not\rightarrow0$ under conditions
(\ref{conds}) and (\ref{condd}), when $\zetaps$ is bounded,
$s_n \zetap\rightarrow0$ is the same as
$s_n \rightarrow0$, which is satisfied if $p = O(n^\lambda)$
with $0 \leq\lambda< 1/2$. When $\zetaps\rightarrow\infty$,
$s_n \zetap\rightarrow0$ is stronger than $s_n \rightarrow0 $.
Under (\ref{conds})--(\ref{condd}),
$\zetap=O(p)$. Hence, the extreme case is $\zetap$
is a constant times $p$, and the condition in part (iv) becomes
$p^2 \sqrt{\log p} / \sqrt{n} \rightarrow0$,
which holds when $p = O(n^\lambda)$ with $0 \leq\lambda< 1/4$.
In the traditional applications with a fixed
$p$, $\zetaps$ is bounded, $s_n \rightarrow
0$ as $n \rightarrow\infty$ and thus
Theorem \ref{thm1} proves that the LDA is
asymptotically optimal.
\end{remark}

The proof of part (iv) of Theorem \ref{thm1} (see Section \ref{sec6}) utilizes
the following lemma, which is also used in the proofs of other results
in this paper.

\begin{lemma}\label{lem1}
Let $\xi_n$ and $\tau_n$ be two sequences of positive numbers
such that $\xi_n \rightarrow\infty$ and $\tau_n \rightarrow0$ as
$n \rightarrow\infty$. If $\lim_{n \rightarrow\infty} \tau_n \xi_n
= \gamma$, where $ \gamma$ may be 0, positive, or $\infty$, then
\[
\lim_{n \rightarrow\infty} \frac{\Phi( - \sqrt{\xi_n }(1-\tau_n))}{
\Phi( - \sqrt{\xi_n} )} = e^{ \gamma} .
\]
\end{lemma}

Since the LDA uses $\bfS^-$ to estimate $\bfSigma^{-1}$ when $p>n$
and is asymptotically
worst as \citet{BicLev04} showed, one may think
that the bad performance of the LDA is caused by the fact that
$\bfS^-$ is not a good estimator of $\bfSigma^{-1}$.
Our following result shows that the LDA may still be
asymptotically worst even if we can estimate
$\bfSigma^{-1}$ perfectly.

\begin{thm}\label{thm2}
Suppose that
(\ref{conds})--(\ref{condd}) hold, $p / n \rightarrow\infty$ and that
$\bfSigma$ is known so that the LDA is given by
(\ref{rule}) with $\hat\bfSigma{}^{-1} = \bfSigma^{-1}$,
$\hat\bfdelta= \bar{\mathbf{x}}_1 - \bar{\mathbf{x}}_2$ and $\hat
{\bar\bfmu
} = \bar{\mathbf{x}}$.
\begin{longlist}[(iii)]
\item[(i)] If $\zetap/\sqrt{p/n} \rightarrow0$ (which is true if
$\zetaps
\not\rightarrow
\infty$), then
$R_{\mathrm{LDA}}(\bfX) \rightarrowp1/2$.
\item[(ii)] If $\zetap/\sqrt{p/n} \rightarrow c$ with $0<c< \infty$,
then $R_{\mathrm{LDA}} (\bfX) \rightarrowp$ a constant strictly
between 0
and $1/2$
and $R_{\mathrm{LDA}} (\bfX)/ \RB\rightarrowp\infty$.
\item[(iii)] If $\zetap/\sqrt{p/n} \rightarrow\infty$, then
$R_{\mathrm{LDA}} (\bfX) \rightarrowp0$ but $R_{\mathrm{LDA}} (\bfX
)/ \RB
\rightarrowp\infty$.
\end{longlist}
\end{thm}

Theorem \ref{thm2} shows that even if $\bfSigma$ is known,
the LDA may be asymptotically
worst and the best we can hope is that the LDA is asymptotically sub-optimal.
It can also be shown that, when $\bfmu_1$ and $\bfmu_2$ are known and
we apply the
LDA with $\hat\bfdelta= \bfdelta$ and $\hat{\bar\bfmu} = ( \bfmu
_1+\bfmu_2)/2$,
the LDA is still not asymptotically optimal when
$\| \bfdelta\|^2 - \| \bfdelta_n \|^2 \not\rightarrow0$,
where $\bfdelta_n$ is any sub-vector of $\bfdelta$ with dimension
$n$.
This indicates that, in order to obtain an asymptotically optimal
classification rule when $p $ is much larger than $n$,
we need sparsity conditions on $\bfSigma$ and $\bfdelta$
when both of them are unknown.
For bounded $\zetaps$ (in which case
the asymptotic optimality is the same as the asymptotic
sub-optimality), by imposing sparsity conditions on
$\bfSigma$, $\bfmu_1$ and $\bfmu_2$, Theorem 2 of
Bickel and Levina (\citeyear{BicLev04}) shows the
existence of an asymptotically optimal classification rule.
In the next section, we obtain a result
by relaxing the boundedness of $\zetaps$ and by imposing sparsity conditions
on $\bfSigma$ and~$\bfdelta$. Since the difference of the two normal
distributions is in $\bfdelta$,
imposing a sparsity condition on $\bfdelta$
is weaker and more reasonable than imposing sparsity
conditions on both $\bfmu_1$ and $\bfmu_2$.

\section{Sparse linear discriminant analysis}\label{sec3}

We focus on the situation where the limit of $p/n$ is positive or
$\infty$.
The following sparsity measure on $\bfSigma$ is considered in \citet{BicLev08}:
%
\begin{equation} \label{cp}
C_{h,p} = \max_{j \leq p} \sum_{l=1}^p |\sigma_{jl}|^h ,
\end{equation}
where $\sigma_{jl}$ is the $(j,l)$th element of $\bfSigma$,
$h$ is a constant not depending on $p$, $0 \leq h < 1$ and $0^0$ is
defined to be 0.
In the special case of $h=0$, $C_{0,p}$
in (\ref{cp}) is the maximum of the numbers of nonzero elements
of rows of $\bfSigma$ so that a $C_{0,p}$ much smaller than
$p$ implies many elements of $\bfSigma$ are equal to 0. If $C_{h,p}$
is much smaller than
$p$ for a constant $h\in(0,1)$,
then $\bfSigma$ is sparse
in the sense that many elements of $\bfSigma$ are very small.
An example of $C_{h,p}$ much smaller than $p$ is $C_{h,p} = O(1)$
or $C_{h,p} = O( \log p )$.

Under conditions (\ref{conds}) and
%
\begin{equation}\label{condp}
\frac{\log p}{n} \rightarrow0,
\end{equation}
\citet{BicLev08} showed that
%
\begin{equation}\label{bl}
\| \tilde{\bfSigma} - \bfSigma\| = O_P ( d_n )
\quad \mbox{and}\quad
\| \tilde\bfSigma{}^{-1} - \bfSigma^{-1} \| = O_P ( d_n
) ,
\end{equation}
where $d_n = C_{h,p} (n^{-1}\log p )^{(1-h)/2}$,
$\tilde{\bfSigma}$ is $\bfS$ thresholded
at $t_n = M_1 \sqrt{\log p}/\sqrt{n}$ with a positive constant $M_1$;
that is, the $(j,l)$th element of $\tilde{\bfSigma}$ is
$\hat\sigma_{jl}I(|\hat\sigma_{jl}| > t_n) $,
$\hat\sigma_{jl}$ is the $(j,l)$th element of $\bfS$ and
$I(A)$ is the indicator function of the set $A$.
We consider a slight modification, that is, only off-diagonal
elements of $\bfS$ are thresholded.
The resulting estimator is still denoted by
$\tilde{\bfSigma}$ and it has property (\ref{bl})
under conditions (\ref{conds}) and (\ref{condp}).

We now turn to the sparsity of $\bfdelta$.
On one hand, a large $\zetaps$ results in a~large difference between
$N_p( \mu_1, \bfSigma)$ and $N_p( \mu_2, \bfSigma)$
so that the optimal rule has a
small misclassification rate. On the other hand,
a larger divergence rate of $\zetaps$
results in a more difficult task of constructing a good
classification rule, since
$\bfdelta$ has to be estimated based on the training sample $\bfX$
of a size much smaller than $p$.
We consider the following sparsity measure on $\bfdelta$ that is
similar to
the sparsity measure $C_{h,p}$ on $\bfSigma$:
%
\begin{equation}
D_{g,p} = \sum_{j=1}^p \delta_{j}^{2g} , \label{dp}
\end{equation}
where $\delta_{j}$ is the $j$th component of $\bfdelta$,
$g$ is a constant not depending on $p$ and $0 \leq g < 1$.
If $D_{g,p}$ is much smaller than\vspace*{-1pt}
$p$ for a $g\in[0,1)$, then $\bfdelta$ is sparse.
For $\zetap$ defined in (\ref{delta}), under (\ref{conds})--(\ref{condd}),
$ \zetap\leq c_0 \|\bfdelta\|^2 \leq c_0^{1+2(1-g)}D_{g,p}$. Hence,
the rate of divergence of $\zetap$ is always smaller than
that of $D_{g,p}$ and, in particular,
$\zetaps$ is bounded when $D_{g,p}$ is bounded for a $g \in[0,1)$.

We consider the sparse estimator $\tilde\bfdelta$ that
is $\hat\bfdelta$ thresholded at
%
\begin{equation}\label{an}
a_n = M_2 \biggl( \frac{\log p}{n} \biggr)^\alpha
\end{equation}
with constants $M_2>0$ and $\alpha\in(0,1/2)$,
that is, the $j$th component of $\tilde\bfdelta$ is
$\hat\delta_jI(|\hat\delta_j | > a_n)$, where
$\hat\delta_j$ is the $j$th component of $\hat\bfdelta$.
The following result is useful.

\begin{lemma}\label{lem2}
Let $\delta_j$ be the $j$th component of $\bfdelta$,
$\hat\delta_j$ be the $j$th component of~$\hat\bfdelta$,
$a_n $ be given by (\ref{an}) and $r >1$ be a fixed constant.
\begin{longlist}[(ii)]
\item[(i)] If (\ref{condp}) holds, then
%
\begin{equation}\label{prob1}
P \biggl( \bigcap_{1 \leq j \leq p,  | \delta_j | \leq a_n/r}\{
| \hat\delta_j | \leq a_n \} \biggr) \rightarrow1
\end{equation}
and
%
\begin{equation}\label{prob2}
P \biggl( \bigcap_{1 \leq j \leq p,  | \delta_j | > ra_n}\{
| \hat\delta_j | > a_n \} \biggr) \rightarrow1 .
\end{equation}
\item[(ii)] Let $q_{n0} = $ the number of $j$'s with $|\delta_j | > r
a_n $,
$q_n = $ the number of $j$'s with $|\delta_j | > a_n /r$ and
$\hat{q} = $ the number of $j$'s with $|\hat\delta_j | > a_n $.
If (\ref{condp}) holds, then
\[
P ( q_{n0} \leq\hat{q} \leq q_n ) \rightarrow1 .
\]
\end{longlist}
\end{lemma}

We propose a sparse linear
discriminant analysis (SLDA) for high-dimen\-sion~$p$,
which is given by (\ref{rule}) with $\hat\delta= \tilde\delta$,
$\hat\bfSigma= \tilde\bfSigma$ and $\hat{\bar\bfmu}=
\bar{\mathbf{x}}$. The following result establishes the asymptotic
optimality of the SLDA under some conditions on the rate of
divergence of $p$, $C_{h,p}$, $D_{g,p}$, $q_n$ and $\zetap$.

\begin{thm}\label{thm3}
Let $C_{h,p}$ be given by (\ref{cp}),
$D_{g,p}$ be given by (\ref{dp}), $a_n $ be given by (\ref{an}),
$q_n$ be as defined in Lemma \ref{lem2} and
$d_n = C_{h,p} (n^{-1}\log p )^{(1-h)/2}$.
Assume that
conditions (\ref{conds}), (\ref{condd}) and (\ref{condp}) hold and
%
\begin{equation}\label{bn}
b_n = \max\Bigl\{ d_n, \frac{a_n^{1-g}\sqrt{D_{g,p} }}{\zetaps}
,  \frac{\sqrt{C_{h,p}q_n}}{\zetaps\sqrt{n}} \Bigr\} \rightarrow0.
\end{equation}
\begin{longlist}[(iii)]
\item[(i)] The conditional misclassification
rate of the SLDA is equal to
\[
R_{\mathrm{SLDA}} (\bfX) =
\Phi\bigl( -[1+O_P(b_n)]\zetaps/ 2 \bigr) .
\]
\item[(ii)] If $\zetaps$ is bounded,
then the SLDA is asymptotically optimal and
\[
\frac{R_{\mathrm{SLDA}} ( \bfX) }{\RB} - 1 = O_P(b_n) .
\]
\item[(iii)] If $\zetaps\rightarrow\infty$, then the SLDA is
asymptotically sub-optimal.
\item[(iv)] If $\zetaps\rightarrow\infty$ and $b_n \zetap
\rightarrow0$, then the SLDA is asymptotically optimal.
\end{longlist}
\end{thm}

\begin{remark}\label{rem2}
Condition (\ref{bn})
may be achieved by an appropriate choice of $\alpha$ in $a_n$, given
the divergence rates of $C_{h,p}$, $D_{g,p}$, $q_n$ and $\zetaps$.
\end{remark}

\begin{remark}\label{rem3}
When $\zetaps$ is bounded and
(\ref{conds})--(\ref{condd}) hold, condition (\ref{bn})
is the same as
%
\begin{equation}\label{cond1}
d_n \rightarrow0,\qquad  D_{g,p} a_n^{2(1-g)} \rightarrow0
\quad \mbox{and}\quad
C_{h,p} q_n /n \rightarrow0 .
\end{equation}
\end{remark}

\begin{remark}\label{rem4}
When $\zetaps\rightarrow\infty$,
condition (\ref{bn}), which is sufficient for the asymptotic
sub-optimality of the SLDA, is implied by
$d_n \rightarrow0$, $D_{g,p} a_n^{2(1-g)} = O(1)$ and
$C_{h,p} q_n /n =O(1)$. When $\zetaps\rightarrow\infty$,
the condition $b_n \zetap\rightarrow0$, which is sufficient
for the asymptotic optimality of the SLDA, is the
same as
%
\begin{equation} \label{cond2}
\zetap d_n \rightarrow0 ,\qquad
\zetap D_{g,p} a_n^{2(1-g)} \rightarrow0 \quad \mbox{and}\quad
\zetap C_{h,p} q_n /n \rightarrow0 .
\end{equation}
\end{remark}

We now study when condition (\ref{bn}) holds and when
$b_n \zetap\rightarrow0$ with $\zetaps\rightarrow\infty$.
By Remarks 3 and 4, (\ref{bn}) is the same as condition (\ref{cond1})
when $\zetaps$ is bounded, and $b_n \zetap\rightarrow0$ is
the same as condition (\ref{cond2}) when $\zetaps\rightarrow\infty$.
\begin{enumerate}
\item
If there are two constants $c_1$ and $c_2$ such that $0 < c_1 \leq
|\delta_j| \leq c_2$
for any nonzero $\delta_j$, then $q_n $ is exactly the number of
nonzero $\delta_j$'s.
Under condition~(\ref{condd}), $\zetap$ and $D_{0,p}$ have exactly
the order $q_n$.
\begin{enumerate}[(a)]
\item[(a)]
If $q_n$ is bounded (e.g., there are only finitely many nonzero $\delta_j$'s), then
$\zetaps$ is bounded and condition (\ref{bn}) is the same as
condition~(\ref{cond1}).
The last two convergence requirements in (\ref{cond1}) are implied
by $d_n = C_{h,p} (n^{-1} \times\log p )^{(1-h)/2} \rightarrow0$, which is
the condition
for the consistency of $\tilde\bfSigma$ proposed by \citet{BicLev08}.
\item[(b)]
When $q_n \rightarrow\infty$ ($\zetaps\rightarrow\infty$), we
assume that $q_n = O(n^\eta)$ and $C_{h,p}= O(n^\gamma)$ with $\eta
\in(0,1)$ and $\gamma\in[0,1)$.
Then, condition (\ref{cond2}) is implied~by
%
\begin{eqnarray} \label{cond22}
n^{\eta+ \gamma} (n^{-1}\log p)^{(1-h)/2} &\rightarrow&0,\qquad
n^{2\eta} (n^{-1}\log p)^{2\alpha}
\rightarrow0,\nonumber\\[-8pt]\\[-8pt]
n^{2\eta+ \gamma-1} &\rightarrow&0.\nonumber
\end{eqnarray}
If we choose $\alpha= (1-h)/4$,
then condition (\ref{cond22}) holds when $2 \eta+ \gamma< 1$
and $n^{\eta+ \gamma} (n^{-1} \log p)^{(1-h)/2} \rightarrow0$.
To achieve (\ref{cond22}) we need to know
the divergence rate of $p$. If $p = O(n^\kappa)$
for a $\kappa\geq1$, then $(n^{-1} \log p)^{(1-h)/2}
= O((n^{-1} \log n )^{(1-h)/2} )$, and thus condition (\ref{cond22})
holds when
$\eta+ \gamma< (1-h)/2$ and $ \eta< (1+h)/2$.
If $p = O(e^{n^\beta})$ for a $\beta\in(0,1)$,
which is referred to as an ultra-high dimension, then
$(n^{-1} \log p)^{(1-h)/2}
= (n^{\beta-1})^{(1-h)/2} $, and
condition (\ref{cond22}) holds if $\eta+ \gamma<(1-h)(1-\beta)/2$
and $\eta< 1- (1-h)(1-\beta)/2$.
\end{enumerate}
\item
Since
\[
\zetap\geq\sum_{j: |\delta_j | > a_n/r} \delta_j^2 \geq q_n
(a_n/r)^2
\]
and
\[
D_{g,p} \geq\sum_{j: |\delta_j | > a_n/r} \delta_j^{2g}
\geq q_n (a_n/r)^{2(1-g)} ,
\]
we conclude that
\begin{equation}\label{qn}
q_n = O \biggl( \min\biggl\{ \frac{\zetap}{a_n^2},
\frac{D_{g,p}}{a_n^{2(1-g)}} \biggr\} \biggr).
\end{equation}
The right-hand side of (\ref{qn}) can be used as a
bound of the divergence rate of $q_n$
when $q_n \rightarrow\infty$, although it may not be a tight bound.
For example,
if $\zetap= O(\log p )$ and
the right-hand side of (\ref{qn}) is used as a bound for $q_n$,
then the last convergence requirement in (\ref{cond1}) or (\ref
{cond2}) is implied
by the first convergence requirement in (\ref{cond1}) or (\ref{cond2})
when $ \alpha\leq(1+h)/4$.
\item
If $D_{g,p} = O( C_{h,p})$, then
the second convergence requirement in (\ref{cond1}) or (\ref{cond2})
is implied by the first convergence requirement in (\ref{cond1}) or
(\ref{cond2})
when $ \alpha\geq(1-h)/[4(1-g)]$.
\item
Consider the case where $C_{h,p}= O(\log p)$,
$D_{g,p}= O(\log p)$ and an ultra-high dimension, that is,
$p = O(e^{n^\beta})$ for a $\beta\in(0,1)$.
From the previous discussion, condition (\ref{cond1}) holds if
$d_n \rightarrow0$, and (\ref{cond2}) holds if
$d_n \log p \rightarrow0$. Since $\log p = O( n^\beta)$,
$d_n = O( n^{\beta+ (\beta-1)(1-h)/2})$, which converges to 0 if
$\beta< (1-h)/(3-h)$. If $\zetaps$ is bounded,
then $d_n \rightarrow0$ is sufficient for condi\-tion~(\ref{bn}).
If $\zetaps\rightarrow\infty$, then the largest divergence rate of
$\zetap$
is $O( \log p) = O( n^\beta)$ and $\zetap d_n \rightarrow0$ (i.e.,
the SLDA is asymptotically optimal) when
$\beta< (1-h)/(5-h)$. When $h =0$, this means $\beta< 1/5$.
\item
If the divergence rate of $p$ is smaller than $O(e^{n^\beta})$
then we can afford to have a larger than $O( \log p)$ divergence
rate for $C_{h,p}$ and $D_{g,p}$.
For~exam\-ple, if $p = O(n^\kappa)$ for a $\kappa\geq1$ and
$\max\{ C_{h,p}, D_{g,p} \} = c n^\gamma$ for a $\gamma\in(0,1)$
and a positive constant $c$, then
$\log p = O( \log n )$ diverges to $\infty$ at a rate slower than
$n^\gamma$. We now study when condition (\ref{cond1}) holds.
First, $d_n = C_{h,p} (n^{-1} \log p )^{(1-h)/2}
= O ( n^{\gamma- (1-h)/2} (\log n )^{(1-h)/2} )$, which converges to 0
if $\gamma< (1-h)/2 \leq1/2$. Second, $a^{2(1-g)} D_{g,p}
= O( n^{\gamma- 2(1-g)\alpha} (\log n )^{2(1-g)\alpha} )$,
which converges to 0 if $\alpha$ is chosen so that $\alpha> \gamma/[2(1-g)]$.
Finally, if we use the right-hand side of (\ref{qn}) as a bound for $q_n$,
then $C_{h,p} q_n /n = O( n^{2 (1-g)\alpha+ \gamma-1} / (\log n
)^{2(1-g)\alpha})$,
which converges to 0 if $\alpha\leq(1-\gamma)/[2(1-g)]$.
Thus, condition (\ref{cond1}) holds if $\gamma< (1-h)/2$ and
$\gamma/[2(1-g)] < \alpha\leq(1-\gamma)/[2(1-g)]$.
For condition (\ref{cond2}), we assume that $\zetap= O(n^{\rho\gamma
})$ with a
$\rho\in[0, 1]$ ($\rho= 0$ corresponds to a bounded $\zetaps$).
Then, a similar analysis leads to the conclusion that condition (\ref{cond2})
holds if $(1+\rho)\gamma\leq(1-h)/2$ and
$(1+\rho) \gamma/[2(1-g)] < \alpha\leq[1-(1+\rho)\gamma]/[2(1-g)]$.
\end{enumerate}

To apply the SLDA, we need to choose two constants, $M_1$ in the thresholding
estimator $\tilde{\bfSigma}$ and $M_2$ in the thresholding estimator
$\tilde\bfdelta$.
We suggest a~data-driven method via a cross-validation procedure.
Let $\bfX_{ki}$ be the data set containing the entire training sample
but with $\mathbf{x}_{ki}$ deleted, and let $T_{ki}$ be the
SLDA rule based on $\bfX_{ki}$, $i=1,\ldots,n_k$, $k=1,2$.
The leave-one-out cross-validation estimator of the misclassification
rate of the SLDA is
\[
\hat{R}_{\mathrm{SLDA}} = \frac{1}{n}
\sum_{k=1}^2 \sum_{i=1}^{n_k} r_{ki} ,
\]
where $r_{ki} $ is the indicator function of whether
$T_{ki}$ classifies $\mathbf{x}_{ki}$ incorrectly.
Let $ R(n_1,n_2)$ denote $R_{\mathrm{SLDA} } $ when the sample sizes are
$n_1$ and $n_2$.
Then
\[
E( \hat{R}_{\mathrm{SLDA}}) = \frac{1}{n}
\sum_{k=1}^2 \sum_{i=1}^{n_k} E( r_{ki} ) =
\frac{n_1R(n_1-1,n_2)+ n_2 R(n_1,n_2-1)}{n} ,
\]
which is close to $R(n_1,n_2)= R_{\mathrm{SLDA}}$ for large $n_k$.
Let $\hat{R}_{\mathrm{SLDA}} (M_1,M_2)$ be the cross-validation
estimator when
$(M_1,M_2)$ is used in thresholding $\hat\bfS$ and $\hat\bfdelta$.
Then, a data-driven method of selecting $(M_1,M_2)$ is to minimize
$\hat{R}_{\mathrm{SLDA}} (M_1,M_2)$ over a suitable range of $(M_1,M_2)$.
The resulting $\hat{R}_{\mathrm{SLDA}} $ can also be used as an estimate
of the misclassification rate of the SLDA.

\section{Extensions}\label{sec4}

We first consider an extension of the main result in Section \ref{sec3}
to nonnormal $\mathbf{x}$ and $\mathbf{x}_{ki}$'s. For nonnormal
$\mathbf{x}$, the
LDA with known $\bfmu_k$ and $\bfSigma$, that is, the
rule classifying $\mathbf{x}$ to class 1 if and only if
$\bfdelta' \bfSigma^{-1} ( \mathbf{x}- \bar{\bfmu}) \geq0$,
is still optimal when $\mathbf{x}$ has an elliptical distribution [see,
e.g., \citet{autokey8}] with density
%
\begin{equation}\label{elli}
c_p | \bfSigma|^{-1/2} f \bigl( (\mathbf{x}- \bfmu)' \bfSigma^{-1} (
\mathbf{x}- \bfmu) \bigr),
\end{equation}
where $\bfmu$ is either $\bfmu_1$ or $\bfmu_2$,
$f$ is a monotone function on $[0, \infty)$,
and $c_p$ is a~normalizing constant.
Special cases of (\ref{elli}) are the multivariate $t$-distribu\-tion
and the multivariate double-exponential distribution.
Although this rule is not necessarily optimal when the distribution of
$\mathbf{x}$ is
not of the form (\ref{elli}), it is still a reasonably good rule when
$\bfmu_k$ and $\bfSigma$
are known. Thus, when $\bfmu_k$ and $\bfSigma$ are unknown,
we study whether the misclassification rate of the SLDA defined in
Section \ref{sec3} is
close to that of the LDA with known $\bfmu_k$ and
$\bfSigma$.\looseness=-1

From the proofs for the asymptotic properties of the SLDA in Section \ref{sec3},
the results depending on the normality assumption are:
\begin{longlist}[(iii)]
\item[(i)] result (\ref{bl}), the consistency of $\tilde\bfSigma$;
\item[(ii)] results (\ref{prob1}) and (\ref{prob2}) in Lemma \ref{lem2};
\item[(iii)] the form of the optimal misclassification rate given by
(\ref
{delta});
\item[(iv)] the result in Lemma \ref{lem1}.
\end{longlist}

Thus, if we relax the normality assumption, we need to address (i)--(iv).
For (i), it was discussed in Section 2.3 of \citet{BicLev08}
that result (\ref{bl}) still holds when the normality assumption is
replaced by one of the following two conditions.
The first condition is
%
\begin{equation}\label{case1}
\sup_{k, j} E( e^{tx_{kij}^2} ) < \infty\qquad \mbox{for all $|t|
\leq t_0$}
\end{equation}
for a constant $t_0 >0$,
where $x_{kij}$ is the $j$th component of $\mathbf{x}_{ki}$. Under
condition~(\ref{case1}),
result (\ref{bl}) holds without any modification. The second condition
is\looseness=-1
%
\begin{equation}\label{case2}
\sup_{k, j} E| x_{kij}|^{2 \nu} < \infty
\end{equation}
for a constant $\nu>0$. Under condition (\ref{case2}), result
(\ref{bl}) holds with $n^{-1} \log p $ changed to $n^{-1} p^{4/\nu}$.
The same argument can be used to address (ii), that is, results
(\ref{prob1}) and (\ref{prob2}) hold under condition (\ref{case1}) or
condition (\ref{case2}) with $n^{-1} \log p $ replaced by $n^{-1}
p^{4/\nu}$.
For (iii), the normality of $\mathbf{x}$ can be relaxed to that, for any
$p$-dimensional nonrandom vector $\bfl$ with $\| \bfl\|=1$ and any
real number $t$,
%
\begin{equation}\label{case3}
P \bigl( \bfl' \bfSigma^{-1/2}( \mathbf{x}- \bfmu) \leq t \bigr) =
\Psi(t ) ,
\end{equation}
where $\Psi$ is an unknown distribution function symmetric about 0 but
it does not depend on $\bfl$. Distributions satisfying (\ref{case3}) include
elliptical distributions [e.g., a distribution of the form (\ref
{elli})] and
the multivariate scale mixture of normals [\citet{autokey8}].
Under (\ref{case3}), when $\bfmu_k$ and $\bfSigma$ are known,
the LDA has misclassification rate
$\Psi( - \zetaps/ 2)$ with $\zetaps$ given by (\ref{delta}).
It remains to address (iv). Note that
the following result,
%
\begin{equation}\label{normal}
\frac{ x }{1+x^2} e^{- x^2/2} \leq\Phi( - x) \leq\frac{1}{x} e^{- x^2/2},\qquad
x >0,
\end{equation}
is the key
for Lemma \ref{lem1}. Without assuming normality, we consider the condition
%
\begin{equation}\label{case4}
0< \lim_{x \rightarrow\infty} \frac{x^\omega e^{- c x^{\varphi}}}{
\Psi( - x ) }
< \infty,
\end{equation}
where $\varphi$ is a constant, $0 \leq\varphi\leq2$,
$\omega$ is a constant and $c$ is a positive constant.
For the case where $\Psi$ is standard normal, condition (\ref{case4})
holds with
$\varphi=2$, $\omega= -1$ and $c =1/2$. Under condition (\ref{case4}),
we can show that the result in Lemma holds for the case of $\gamma=
0$, which
is needed to extend the result in Theorem \ref{thm3}(iv).
This leads to the following extension.

\begin{thm}\label{thm4}
Assume condition (\ref{case3}) and either
condition (\ref{case1})
or (\ref{case2}). When condition (\ref{case1}) holds, let $b_n$ be
defined by (\ref{bn}).
When condition (\ref{case2}) holds,
let $a_n$ and $b_n$ be defined by (\ref{an}) and (\ref{bn}),
respectively, with $n^{-1}\log p$ replaced by $n^{-1} p^{4/\nu}$.
Assume that $a_n \rightarrow0$ and $b_n \rightarrow0$.
\begin{longlist}[(iii)]
\item[(i)] The conditional misclassification rate of the SLDA is
\[
R_{\mathrm{SLDA}} ( \bfX) = \Psi\bigl( -[1+O_P(b_n)] \zetaps/2 \bigr).
\]
\item[(ii)] If $\zetaps$ is bounded, then
\[
\frac{R_{\mathrm{SLDA}}(\bfX) }{\Psi( - \zetaps/2 )} - 1 = O_P(
b_n) ,
\]
where $\Psi( - \zetaps/2 )$ is the misclassification rate of the LDA
when $\bfmu_k$ and $\bfSigma$ are known.
\item[(iii)] If $\zetaps\rightarrow\infty$, then
$ R_{\mathrm{SLDA}}(\bfX) \rightarrowp0$.
\item[(iv)] If $\zetaps\rightarrow\infty$ and $b_n \zetap
\rightarrow0$,
then
\[
\frac{R_{\mathrm{SLDA}}(\bfX) }{\Psi( - \zetaps/2 )} \rightarrowp
1 .
\]
\end{longlist}
\end{thm}

We next consider extending the results in Sections \ref{sec2} and \ref
{sec3} to
the classification problem with $K \geq3$ classes.
Let $\mathbf{x}$ be a $p$-dimensional normal random vector belonging to
class $k$ if $\mathbf{x}\sim N_p( \bfmu_k , \bfSigma)$, $k=1,\ldots,K$,
and the training sample be $\bfX= \{ \mathbf{x}_{ki}, i=1,\ldots,n_k,
k=1,\ldots,K \}$, where
$n_k$ is the sample size for class $k$,
$\mathbf{x}_{ki} \sim N_p( \bfmu_k , \bfSigma)$, $k=1,\ldots,K$,
and all $\mathbf{x}_{ki}$'s are independent.
The LDA classifies $\mathbf{x}$ to class $k$ if and only if
$\hat{\bfdelta}{}_{kl}' \hat\bfSigma{}^{-1} ( \mathbf{x}- \hat{\bar
{\bfmu}}_{kl} ) \geq0$
for all $l \neq k$, $l=1,\ldots,K$, where
$\hat{\bfdelta}_{kl} = \bar{\mathbf{x}}_k - \bar{\mathbf{x}}_l$,
$\hat{\bar{\bfmu}}_{kl} = (\bar{\mathbf{x}}_k+\bar{\mathbf{x}}_l)/2$,
$\bar{\mathbf{x}}_k = n_k^{-1} \sum_{i=1}^{n_k} \mathbf{x}_{ki}$ and
$\hat\bfSigma{}^{-1} $ is an inverse or a generalized inverse of $\bfS
= n^{-1} \sum_{k=1}^K
\sum_{i=1}^{n_k} ( \mathbf{x}_{ki} - \bar{\mathbf{x}}_k )( \mathbf
{x}_{ki} -
\bar{\mathbf{x}}_k )'$,
and $n=n_1+\cdots+ n_K$. The conditional misclassification rate of the
LDA is
\[
\frac{1}{K} \sum_{k=1}^K \sum_{j \neq k} P_k \bigl(
\hat{\bfdelta}{}_{jl}' \hat\bfSigma{}^{-1} ( \mathbf{x}- \hat{\bar
{\bfmu}}_{jl} ) \geq0 , l \neq j \bigr) ,
\]
where $P_k$ is the probability with respect to $\mathbf{x}\sim N_p(
\bfmu
_k , \bfSigma)$,
$k=1,\ldots,K$. The SLDA and its conditional misclassification rate
can be obtained by simply replacing $\hat{\bfSigma}$ and $\hat
{\bfdelta}_{kl}$
by their thresholding estimators $\tilde{\bfSigma}$ and $\tilde
{\bfdelta}_{kl}$,
respectively. For simplicity of computation, we suggest the use of
the same thresholding constant (\ref{an}) for all $\tilde{\bfdelta}_{kl}$'s.

The optimal rate can be calculated as
%
\begin{equation}\label{misrate}
\RB= \frac{1}{K} \sum_{k=1}^K \sum_{j \neq k} P_k \bigl(
\bfdelta_{jl}' \bfSigma^{-1} ( \mathbf{x}- \bar{\bfmu}_{jl} ) \geq
0 ,
l \neq j
\bigr) ,
\end{equation}
where $\bfdelta_{jl} = \bfmu_j - \bfmu_l$ and
$\bar{\bfmu}_{jl} = (\bfmu_j+\bfmu_l)/2$, $j,l =1,\ldots,K$, $j
\neq l$.
Asymptotic properties of the LDA and SLDA can be obtained,
under the asymptotic setting with $n \rightarrow\infty$ and
$ n_k / n \rightarrow$ a constant in $(0,1)$ for each $k$.
Sparsity conditions should be imposed to each $\bfdelta_{kl}$.
If the probabilities in expression~(\ref{misrate}) do not converge to
0, then
the asymptotic optimality of the LDA (under the conditions in Theorem
\ref{thm1})
or the SLDA (under the conditions in Theorem \ref{thm3}) can be established
using the
same proofs as those in Section \ref{sec6}. When
$\RB$ in (\ref{misrate}) converges to 0,
to consider convergence rates,
the proof of the asymptotic optimality of the LDA or SLDA requires an
extension of Lemma \ref{lem1}.
Specifically, we need an extension of result (\ref{normal})
to the case of multivariate normal distributions. This
technical issue, together with empirical properties of the SLDA
with $K \geq3$, will be investigated in our future research.

\section{Numerical studies}\label{sec5}

Golub et al. (\citeyear{Goletal99}) applied gene expression microarray techniques to study
human acute leukemia and discovered the distinction between acute myeloid
leukemia (AML) and acute lymphoblastic leukemia (ALL). Distinguishing ALL
from AML is crucial for successful treatment, since chemotherapy regimens
for ALL can be harmful for AML patients.
An accurate classification based solely on gene expression
monitoring independent of previous biological knowledge is desired
as a general strategy for discovering and predicting cancer classes.

We considered a dataset that was used by many researchers
[see, e.g., \citet{FanFan08}]. It
contains the expression levels of $p=7\mbox{,}129$ genes for
$n=72$ patients.
Patients in the sample are known to come from two
distinct classes of leukemia: $n_1=47$ are from the ALL class, and
$n_2=25$ are from the AML class.

\begin{figure}

\includegraphics{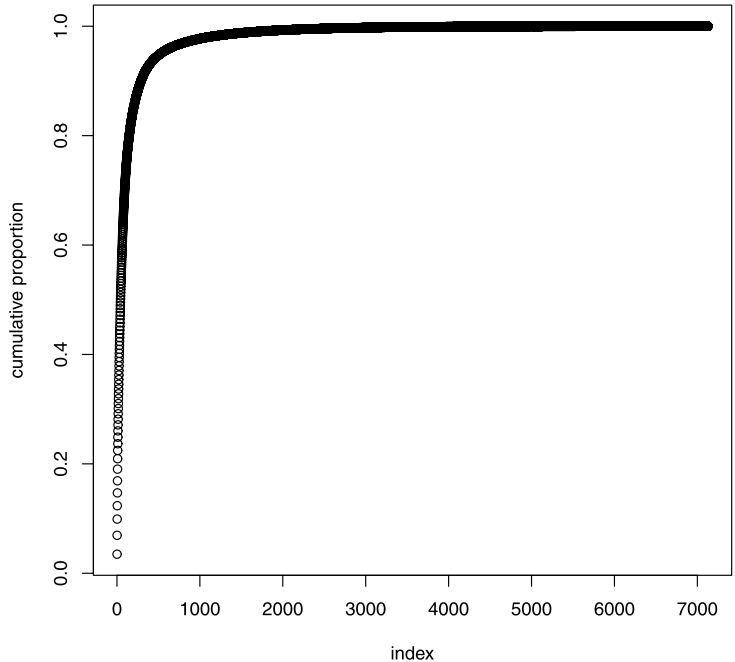}

\caption{Cumulative proportions.}\label{fig1}
\end{figure}

\begin{figure}

\includegraphics{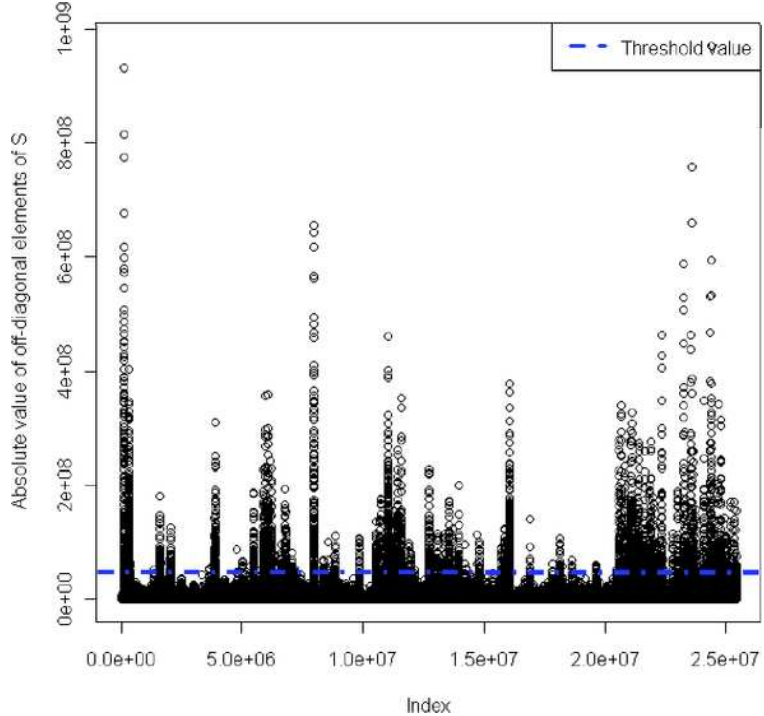}

\caption{Plot of off-diagonal elements of $\bfS$.}\label{fig2}
\end{figure}

Figure \ref{fig1} displays the cumulative proportions defined as\vspace*{-2pt}
$\sum_{j=1}^l \hat\delta_{(j)}^2 / \| \hat\bfdelta\|^2$,
$l=1,\ldots,p$,
where $\hat\delta_{(j)}^2$ is the $j$th largest value among the
squared components\vspace*{-2pt}
of $\hat\bfdelta$.
These proportions indicate the importance of the contribution of each
$\hat\delta_{(j)}$.
It can be seen from Figure \ref{fig1} that the first 1,000 $\hat\delta_{(j)}$'s
contribute a~cumulative proportion nearly 98\%. Figure \ref{fig2} plots the absolute values of
the off-diagonal elements of
the sample covariance matrix $\bfS$. It can be seen that many of
them are relatively small. If we ignore a factor of $10^8$, then among
a total of 25,407,756 values in Figure \ref{fig2}, only 0.45\%
of them vary from 0.35 to 9.7
and the rest of them are under 0.35.

\begin{figure}

\includegraphics{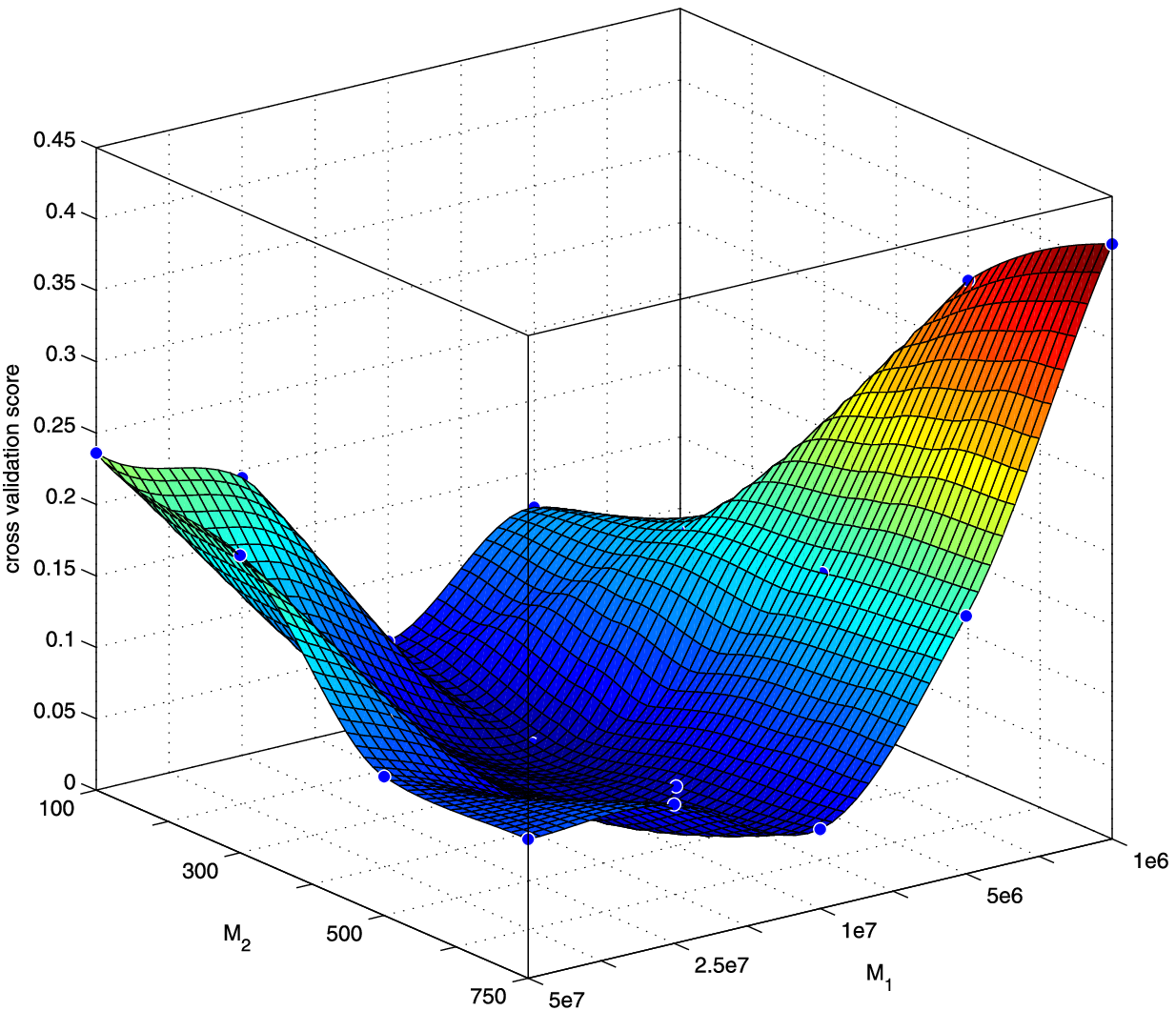}

\caption{Cross-validation score vs $(M_1,M_2)$.}\label{fig3}
\end{figure}

For the SLDA, to construct sparse estimates of $\bfdelta$ and
$\bfSigma$ by thresholding,
we applied the cross-validation method described in the end of Section
\ref{sec3} to choose the
constants $M_1$ and $M_2$ in the thresholding values $t_n
= M_1 (n^{-1} \log p)^{0.5}$ and $a_n = M_2 (n^{-1} \log p )^{0.3}$.
Figure \ref{fig3} shows the cross validation scores
$\hat{R}_{\mathrm{SLDA}} (M_1,M_2)$ over a range of $(M_1,M_2)$.
The minimum cross validation score is achieved at $M_1 = 10^7$ and
$M_2 = 300$. These thresholding values resulted in a $\tilde\bfdelta$ with
exactly 2,492 nonzero components, which is about 35\% of all
components of $\hat\bfdelta$,
and a $\tilde\bfSigma$ with exactly 227,083
nonzero elements, which is about 0.45\% of all elements of $\bfS$.
Note that the number of nonzero estimates of $\bfdelta$ is still much larger
than $n=72$, but the SLDA does not require it to be smaller than $n$.
The resulting SLDA has an estimated (by cross validation)
misclassification rate
0.0278. In fact, 1 of the 47 ALL cases and 1 of the 25 AML cases
are misclassified under the cross validation evaluation of the SLDA.

For comparison, we carried out the LDA with a generalized inverse
$\bfS^-$. In the leave-one-out cross-validation evaluation of the LDA,
2 of the 47 ALL cases and 5 of the 25 AML cases
are misclassified by the LDA, which results in an estimated
misclassification rate 0.0972. Compared with the LDA,
the SLDA reduces the misclassification rate
by nearly 70\%. From Figure 5 of \citet{FanFan08}, the
misclassification rate of the FAIR method, estimated by the average of 100
randomly constructed
cross validations with $\pi n$ data points for constructing classifier and
$(1-\pi)n$ data points for validation ($\pi= 0.4, 0.5$ and $0.6$),
ranges from 5\% to 7\%, which is smaller than
the misclassification rate of the LDA but larger than
the misclassification rate of the SLDA.

\begin{figure}

\includegraphics{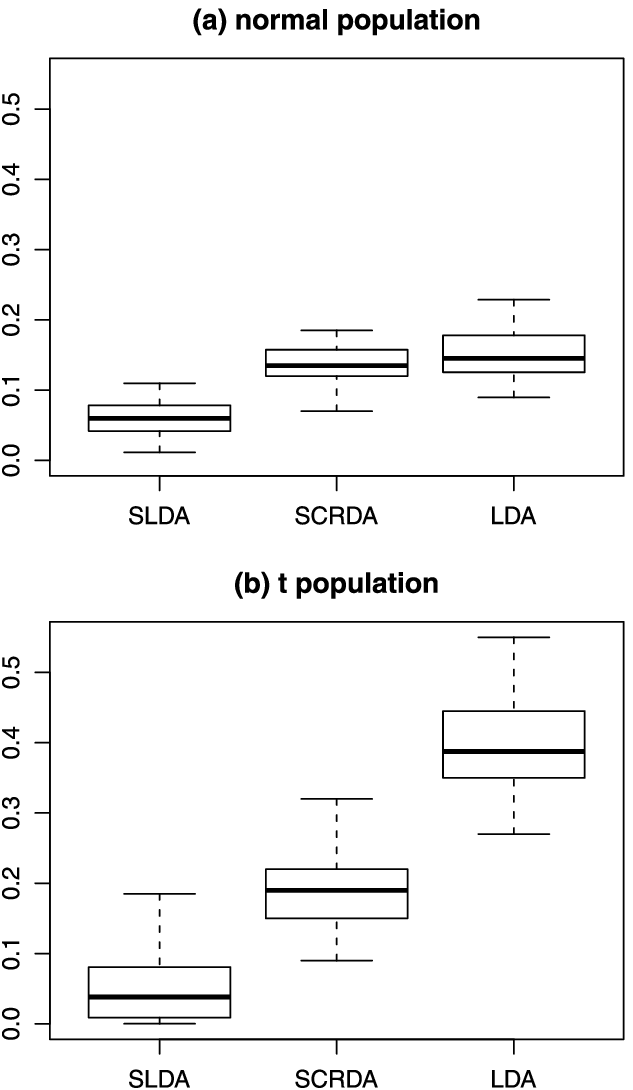}

\caption{Boxplots of conditional misclassification rates of SLDA, SCRDA
and LDA.}\label{fig4}
\end{figure}

We also performed a simulation
study on the conditional misclassification rate of SLDA under a
population constructed using estimates from the real data set and a
smaller dimension $p=1\mbox{,}714$. The smaller dimension was used
to reduce the computational cost and the
1,714 variables were chosen
from the 7,129 variables with
$p$-values (of the two sample $t$-tests for the mean effects)
smaller than 0.05.
In each of the 100 independently generated data sets,
independent $\{ \mathbf{x}_{1i}, i=1,\ldots,47\}$ and $\{ \mathbf{x}_{2i},
i=1,\ldots,25\}$ were generated from $N_p ( \hat\bfmu_1, \tilde
\bfSigma
)$ and $N_p ( \hat\bfmu_2, \tilde\bfSigma)$, respectively, where
$p=1\mbox{,}714$ and $\hat\bfmu_k$ and $\tilde\bfSigma$ are estimates from
the real data set. The sparse estimate $\tilde\bfSigma$ was used
instead of the sample covariance matrix $\bfS$, because $\bfS$ is
not positive definite.
Since the population means and covariance
matrix are known in the simulation, we were able to compute the
conditional misclassification rate $R_{\mathrm{SLDA}} ( \bfX)$
for each generated data set.
A boxplot of 100 values of $R_{\mathrm{SLDA}} (
\bfX)$ in the simulation is given in Figure \ref{fig4}(a).
The unconditional misclassification rate of the
SLDA can be approximated by
averaging over the 100 conditional misclassification rates. In this
simulation, the unconditional misclassification rate for the SLDA
is 0.069. Since the population is known in simulation,
the optimal misclassification rate $R_{\mathrm{OPT}}$ is known to be 0.03.

For comparison, in the simulation we computed the conditional
misclassification rates, $R_{\mathrm{LDA}} (\bfX)$ for the LDA and
$R_{\mathrm{SCRDA}} ( \bfX)$ for the shrunken centroids regularized
discriminant analysis (SCRDA)
proposed by \citet{GuoHasTib07}. Since $R_{\mathrm{SCRDA}} ( \bfX)$
does not have an explicit form, it is approximated
by an independent test data set of size $100$ in each simulation run.
Boxplots of $R_{\mathrm{LDA}} ( \bfX)$ and
$R_{\mathrm{SCRDA}} ( \bfX)$ for 100 simulated data sets are included
in Figure \ref{fig4}(a). It can be seen that the conditional
misclassification rate of the LDA varies more than that of the SLDA.
The unconditional misclassification rate for the LDA,
approximated by the 100
simulated $R_{\mathrm{LDA}} (\bfX)$ values, is 0.152, which
indicates a 53\% improvement of the SLDA over
the LDA in terms of the unconditional
misclassification rate. The SCRDA has a simulated
unconditional misclassification rate 0.137 and its performance is better
than that of the LDA but worse than that of the SLDA.
In this simulation, we also found that
the conditional misclassification rate of the FAIR method
was similar to that of the LDA.

To examine the performance of these classification methods in the
case of nonnormal data, we repeated the same simulation with
the multivariate normal distribution replaced by the multivariate
$t$-distribution with 3 degrees of freedom. The boxplots are given in
Figure \ref{fig4}(b) and the simulated
unconditional misclassification rates are 0.059, 0.194 and 0.399
for the SLDA, SCRDA and LDA, respectively.
Since the $t$-distribution has a larger variability than the
normal distribution, all conditional misclassification rates in the
$t$-distribution case
vary more than those in the normal distribution case.\looseness=-1

\section{Proofs}\label{sec6}
\mbox{}

\begin{pf*}{Proof of Theorem \ref{thm1}}
(i) Let $\hat\sigma_{j,l}$ and $\sigma
_{j,l}$ be
the $(j,l)$th elements of $\bfS$ and $\bfSigma$, respectively.
From result (10) in \citet{BicLev08},
$ \max_{j, l \leq p} | \hat\sigma_{j,l} - \sigma_{j,l} | = O_P
( \sqrt{\log p}
/ \sqrt{n} )$.
Then,
\[
\| \bfS- \bfSigma\| \leq\max_{j \leq p} \sum_{l=1}^p
| \hat\sigma_{j,l} - \sigma_{j,l} | = O_P \bigl( p \sqrt{\log p}
/ \sqrt{n}\, \bigr) = O_P( s_n),
\]
where $\| \bfA\|$ is the norm of the matrix $\bfA$ defined as the
maximum of all eigenvalues of $\bfA$. By (\ref{conds})--(\ref
{condd}) and
$s_n \rightarrow0$, $\bfS^{-1}$ exists and
\[
\| \bfS^{-1} - \bfSigma^{-1} \| =
\| \bfS^{-1} ( \bfS- \bfSigma) \bfSigma^{-1} \| \leq
\| \bfS^{-1} \| \| \bfS- \bfSigma\| \| \bfSigma^{-1} \|
= O_P ( s_n ) .
\]
Consequently,
\[
\hat\bfdelta{}' \bfS^{-1} \bfSigma\bfS^{-1} \hat\bfdelta
= \hat\bfdelta{}' \bfS^{-1} \hat\bfdelta[1+O_P(s_n)]
= \hat\bfdelta{}' \bfSigma^{-1} \hat\bfdelta[1+O_P(s_n)] .
\]
Since $E [( \hat\bfdelta- \bfdelta)' \bfSigma^{-1} ( \hat\bfdelta
- \bfdelta)]
= O(p/n)$ and\vspace*{-1pt} $E[\bfdelta' \bfSigma^{-1} ( \hat\bfdelta- \bfdelta)]^2
\leq\zetap E[( \hat\bfdelta- \bfdelta)' \bfSigma^{-1} (
\hat\bfdelta- \bfdelta)]$,
we have
\begin{eqnarray*}
\hat\bfdelta{}' \bfSigma^{-1} \hat\bfdelta& = &
\bfdelta' \bfSigma^{-1} \bfdelta+ 2 \bfdelta' \bfSigma^{-1} ( \hat
\bfdelta- \bfdelta)
+( \hat\bfdelta- \bfdelta)' \bfSigma^{-1} ( \hat\bfdelta-
\bfdelta)\\
& = & \zetap+
O_P \biggl( \frac{\sqrt{ p}\zetaps}{\sqrt{n}} \biggr) + O_P \biggl(
\frac{p}{n} \biggr)\\
& = & \zetap\biggl[1 + O_P \biggl( \frac{\sqrt{p}}{\sqrt
{n}\zetaps}
\biggr) + O_P \biggl( \frac{p}{n \zetap} \biggr) \biggr]\\
& = & \zetap[1 + O_P ( s_n ) ],
\end{eqnarray*}
where the last equality follows from
$\sqrt{p}/(s_n \sqrt{n} \zetaps) =
1 / (\sqrt{ p \log p }\zetaps) = O(1)$.
Combining these results, we obtain that
\begin{eqnarray*}
\hat\bfdelta\bfS^{-1} \hat\bfdelta&=&
\hat\bfdelta{}' \bfSigma^{-1} \hat\bfdelta[1+O_P(s_n)]
= \zetap[1+O_P(s_n)]^2\\
& =& \zetap[1+O_P(s_n)].
\end{eqnarray*}
Then
\begin{eqnarray*}
\frac{\hat\bfdelta{}' \bfS^{-1} ( \bar{\mathbf{x}}_1-\bfmu_1 )
- \hat\bfdelta{}' \bfS^{-1} \hat\bfdelta/2}{\sqrt{\hat\bfdelta{}'
\bfS^{-1} \bfSigma\bfS^{-1} \hat\bfdelta}} & = &
- \frac{\sqrt{\hat\bfdelta{}'\bfS^{-1} \hat\bfdelta}}{2\sqrt{1+O_P(s_n)}}
+ \frac{\hat\bfdelta{}' \bfS^{-1} ( \bar{\mathbf{x}}_1-\bfmu_1 ) }{
\sqrt{\hat\bfdelta{}'\bfS^{-1} \bfSigma\bfS^{-1} \hat\bfdelta}}\\
& = & - \frac{\sqrt{\zetap[1+O_P(s_n)]}}{2\sqrt{1+O_P(s_n)}}
+ O_P \biggl( \sqrt{\frac{p}{n}} \biggr)\\
& = & - \frac{\zetaps}{2} [ 1+ O_P(s_n) ]+
O_P \biggl( \sqrt{\frac{p}{n}} \biggr)\\
& = & - \frac{\zetaps}{2} \biggl[ 1+ O_P(s_n) + O_P\biggl (
\frac{\sqrt{p}}{\sqrt{n} \zetaps} \biggr) \biggr]\\
& = & - \frac{\zetaps}{2} [ 1+ O_P(s_n) ].
\end{eqnarray*}
Similarly, we can show that
\[
\frac{\hat\bfdelta{}' \bfS^{-1} ( \bfmu_2 - \bar{\mathbf{x}}_2 )
- \hat\bfdelta{}' \bfS^{-1} \hat\bfdelta/2}{\sqrt{\hat\bfdelta{}'
\bfS^{-1} \bfSigma\bfS^{-1} \hat\bfdelta}} = - \frac{\zetaps}{2}
[ 1+ O_P(s_n) ].
\]
These results and formula (\ref{rate}) imply the result in (i).
{\smallskipamount=0pt
\begin{longlist}[(iii)]
\item[(ii)] Let $\phi$ be the density of $\Phi$. By the result in (i),
\[
R_{\mathrm{LDA}}(\bfX) - \RB= \phi( \omega_n ) O_P(s_n) ,
\]
where $\omega_n$ is between $- \zetaps/2$ and $- [1+O_P(s_n)] \zetaps/2$.
Since $\phi(\omega_n ) $ is bounded by a constant,
the result follows from the fact that $\RB$ is bounded away from 0
when $\zetaps$ is bounded.

\item[(iii)] When $\zetaps\rightarrow\infty$, $\RB\rightarrow0$, and, by
the result in (i),
$R_{\mathrm{LDA}}(\bfX) \rightarrowp0$.

\item[(iv)] If $\zetaps\rightarrow\infty$, then, by
Lemma \ref{lem1} and the condition $ s_n \zetap\rightarrow0$, we
conclude that
$R_{\mathrm{LDA}}(\bfX)/\RB\rightarrowp1$.\qed
\end{longlist}}
\noqed
\end{pf*}

\begin{pf*}{Proof of Lemma \ref{lem1}}
It follows from result (\ref{normal}) that
\begin{eqnarray*}
\frac{\xi_n (1 - \tau_n)}{1+\xi_n (1-\tau_n)^2}
e^{[\xi_n - \xi_n(1-\tau_n)^2]/2} & \leq&
\frac{\Phi( - \sqrt{\xi_n} (1-\tau_n))}{
\Phi( - \sqrt{\xi_n} )} \\
& \leq& \frac{1+ \xi_n}{\xi_n (1 - \tau_n)}
e^{[\xi_n - \xi_n(1-\tau_n)^2]/2} .
\end{eqnarray*}
Since $\xi_n \rightarrow\infty$ and $\tau_n \rightarrow0$,
\[
\frac{\xi_n (1 - \tau_n)}{1+\xi_n (1-\tau_n)^2} \rightarrow1
\quad\mbox{and}\quad \frac{1+ \xi_n}{\xi_n (1 - \tau
_n)}\rightarrow1.
\]
The result follows from $[\xi_n - \xi_n(1-\tau_n)^2]/2 =
\xi_n \tau_n (1 - \tau_n/2) \rightarrow\gamma$ regardless of whether
$\gamma$ is 0, positive, or $\infty$.
\end{pf*}

\begin{pf*}{Proof of Theorem \ref{thm2}}
For simplicity, we prove the case of
$n_1=n_2=n/2$.
\begin{longlist}[(iii)]
\item[(i)] The conditional
misclassification rate of the LDA in this case is given by
(\ref{rate}) with $\hat{\bfSigma}$ replaced by $\bfSigma$.
Note that $\bfSigma^{-1/2} ( \bar{\mathbf{x}}_k - \bfmu_k ) \sim N_p(
\bfzero, n_1^{-1} \bfI)$,
where $\bfI$ is the identity matrix of order $p$.
Let $\zeta_j$ be the $j$th component of $\bfSigma^{-1/2} \bfdelta$.
Then, $\sum_{j=1}^p \zeta_j^2 = \zetap$ and
the $j$th component of
$\bfSigma^{-1/2} ( \bar{\mathbf{x}}_k - \bfmu_k ) $ is $n^{-1/2}_1
\varep
_{kj}$, and
the $j$th component of $\bfSigma^{-1/2} \hat\bfdelta$ is $\zeta_j +
n_1^{-1/2}
(\varep_{1j} - \varep_{2j})$, $j=1,\ldots,p$, where
$\varep_{kj}$, $j=1,\ldots,p$, $k=1,2$, are independent standard normal
random variables.
Consequently,
\begin{eqnarray*}
\hat\bfdelta{}' \bfSigma^{-1} ( \bar{\mathbf{x}}_1-\bfmu_1 )
- \hat\bfdelta{}' \bfSigma^{-1} \hat\bfdelta/2 & = &
\sum_{j=1}^p \biggl( - \frac{\zeta_j^2}{2}
+ \frac{\varep_{1j}^2 - \varep_{2j}^2}{n}
+ \frac{\zeta_j \varep_{2j}}{\sqrt{n_1}} \biggr)\\
& = & - \frac{\zetap}{2}
+ \frac{1}{n} \sum_{j=1}^p (\varep_{1j}^2 - \varep_{2j}^2)
+ \frac{1}{\sqrt{n_1}} \sum_{j=1}^p \zeta_j \varep_{2j}\\
& = & - \frac{\zetap}{2}
+ O_P \biggl( \frac{\sqrt{p}}{n} \biggr) +
O_P \biggl( \frac{ \zetaps}{\sqrt{n}} \biggr)
\end{eqnarray*}
and
\begin{eqnarray*}
\hat\bfdelta{}' \bfSigma^{-1} \hat\bfdelta& = &
\sum_{j=1}^p \biggl( \zeta_j + \frac{ \varep_{1j} - \varep
_{2j}}{\sqrt{n_1}}
\biggr)^2 \\
& = & \zetap
+ \frac{1}{n_1} \sum_{j=1}^p ( \varep_{1j}-\varep_{2j})^2
+ \frac{2}{\sqrt{n_1}} \sum_{j=1}^p \zeta_j ( \varep_{1j}-\varep
_{2j})\\
& = & \zetap
+ \frac{4p}{n} [1+ o_P(1)] +
O_P \biggl( \frac{ \zetaps}{\sqrt{n}} \biggr)\\
& = & \zetap+ \frac{4p}{n} [1+ o_P(1)] ,
\end{eqnarray*}
where the last equality follows from $\zetap= O(p)$ under
(\ref{conds})--(\ref{condd}). Combining these results, we obtain that
%
\begin{equation}\label{r1}
\qquad \frac{ \hat\bfdelta{}' \bfSigma^{-1} ( \bar{\mathbf{x}}_1-\bfmu_1 )
- \hat\bfdelta{}' \bfSigma^{-1} \hat\bfdelta/2}{\sqrt{\hat\bfdelta{}
' \bfSigma^{-1} \hat\bfdelta}}
= - \frac{\zetap}{2
\sqrt{ \zetap+ ({4p}/{n}) [1+ o_P(1)] } } + o_P(1).
\end{equation}
Similarly, we can prove that (\ref{r1}) still holds if
$\bar{\mathbf{x}}_1 - \bfmu_1$ is replaced by $\bfmu_2 - \bar
{\mathbf{x}}_2$.
If $\zetap/\sqrt{p/n} \rightarrow0$, then the quantity in (\ref{r1})
converges to 0 in probability. Hence, $R_{\mathrm{LDA}} (\bfX
)\rightarrowp
1/2$.

\item[(ii)] Since $p/n \rightarrow\infty$, $\zetap/ (p/n)
\rightarrow0$.
Then, the quantity in (\ref{r1}) converges to
$- c/4$ in probability and, hence,
$R_{\mathrm{LDA}} (\bfX) \rightarrowp\Phi(- c/4 ) $,
which is a constant between 0 and $1/2$. Since $\zetaps
\rightarrow\infty$, $\RB\rightarrow0$ and, hence,
$R_{\mathrm{LDA}} (\bfX)/ \RB\rightarrowp\infty$.

\item[(iii)] When $\zetap/\sqrt{p/n} \rightarrow\infty$,
it follows from (\ref{r1}) that the
quantity on the left-hand side of (\ref{r1}) diverges
to $- \infty$ in probability. This proves that $R_{\mathrm{LDA}}
(\bfX)
\rightarrowp0$.
To show $R_{\mathrm{LDA}} (\bfX)/ \RB$ $ \rightarrowp\infty$, we
need a
more refined analysis.
The quantity on the left-hand side of (\ref{r1}) is equal to
\[
- \frac{\zetap+ O_P (\sqrt{p}/n) + O_P( \zetaps/\sqrt{n}) }{2
\sqrt{ \zetap+ ({4p}/{n}) [1+ o_P(1)] } } =
- \frac{\zetaps}{2} (1-\tau_n ),
\]
where
\[
\tau_n = 1- \frac{ \zetaps
+ O_P( \sqrt{p}/n) /\zetaps
+ O_P( 1/\sqrt{n} ) }{\sqrt{\zetap+ ({4p}/{n}) [1+ o_P(1)] } }
\]
and $P( 0 \leq\tau_n \leq1 )\rightarrow1$. Note that
\begin{eqnarray*}
\tau_{1n} & = &
1- \frac{\zetaps}{\sqrt{\zetap+ ({4p}/{n}) [1+ o_P(1)] } }\\
& = &
\frac{({4p}/{n}) [1+ o_P(1)] }{\zetap+ ({4p}/{n}) [1+ o_P(1)] +
\zetaps\sqrt{\zetap+ ({4p}/{n}) [1+ o_P(1)] } }
\end{eqnarray*}
and
\begin{eqnarray*}
\tau_{2n} &=&
\frac{ O_P( \sqrt{p}/n) /\zetaps+ O_P( 1/\sqrt{n} ) }{
\sqrt{\zetap+ ({4p}/{n}) [1+ o_P(1)] } }
= \frac{ O_P( \sqrt{p}/n)}{\zetap} + \frac{O_P( 1/\sqrt{n} )
}{\zetaps}\\
&=& \frac{O_P( \sqrt{p/n} )}{\zetap}
\end{eqnarray*}
under (\ref{conds}) and (\ref{condd}). Then
\[
\tau_n \zetap
= \tau_{1n} \zetap+ \tau_{2n} \zetap
= \tau_{1n} \zetap+ O_P\bigl( \sqrt{p/n} \bigr) .
\]
If $\zetap/ (p/n)$ is bounded, then
$\tau_{1n} \geq c$ for a constant $c>0$ and
\[
\tau_{n} \zetap\geq c \zetap+ O_P\bigl( \sqrt{p/n} \bigr) ,
\]
which diverges to $\infty$ in probability since
$\zetap/ \sqrt{p/n} \rightarrow\infty$.
If $\zetap/ (p/n)\rightarrow\infty$, then
$\tau_{1n}\zetap\geq c p/n$ for a constant $c>0$ and
\[
\tau_{n} \zetap\geq c p/n + O_P\bigl( \sqrt{p/n} \bigr) ,
\]
which diverges to $\infty$ in probability since $p/n \rightarrow
\infty$.
Thus, $\tau_{n}\zetap\rightarrow\infty$ in probability, and
the result follows from Lemma \ref{lem1}.\qed
\end{longlist}
\noqed
\end{pf*}

\begin{pf*}{Proof of Lemma \ref{lem2}}
(i) It follows from (\ref{normal}) that, for
all $t$,
\[
P( | \hat\delta_j - \delta_j | > t ) \leq c_1 e^{-c_2 n t^2} ,
\]
where $c_1$ and $c_2$ are positive constants. Then, the probability in
(\ref{prob1}) is
\begin{eqnarray*}
1- P \biggl( \bigcup_{1 \leq j \leq p,  | \delta_j | \leq a_n/r}
\{ | \hat\delta_j | > a_n \} \biggr) & \geq&
1 - \sum_{j=1}^p P \bigl( | \hat\delta_j - \delta_j | > a_n(r-1)/r
\bigr)\\
& \geq& 1 - p c_1 e^{-c_2 n a_n^2(r-1)^2/r^2 } .
\end{eqnarray*}
Because
\[
\frac{n a_n^2}{\log p} = \biggl( \frac{n}{\log p }
\biggr)^{1-2\alpha}
\rightarrow\infty
\]
when $\alpha< 1/2$, we conclude that $p c_1 e^{-c_2 n a_n^2(r-1)^2/r^2 }
\rightarrow0$, and thus (\ref{prob1}) holds.
The proof of (\ref{prob2}) is similar since
\begin{eqnarray*}
1- P \biggl( \bigcup_{1 \leq j \leq p,  | \delta_j | > ra_n}\{
| \hat\delta_j | \leq a_n \} \biggr) & \geq&
1 - \sum_{j=1}^p P \bigl( | \hat\delta_j - \delta_j | > a_n(r-1)
\bigr)\\
& \geq& 1 - p c_1 e^{-c_2 n a_n^2 (r-1)^2} .
\end{eqnarray*}
(ii) The result follows from results (\ref{prob1}) and (\ref{prob2}).
\end{pf*}

\begin{pf*}{Proof of Theorem \ref{thm3}}
The conditional misclassification rate
$R_{\mathrm{SLDA}}(\bfX)$ is
given by
\[
\frac{1}{2} \sum_{k=1}^2 \Phi \biggl(
\frac{(-1)^k \tilde\bfdelta {}' \tilde{\bfSigma}^{-1} ( \bfmu_k -\bar{\bfx}_k )
- \hat\bfdelta{}' \tilde{\bfSigma}^{-1} \tilde\bfdelta /2}{\sqrt{\tilde\bfdelta {}'
\tilde{\bfSigma}^{-1} \bfSigma \tilde{\bfSigma}^{-1} \tilde\bfdelta }} \biggr).
\]
From result (\ref{bl}),
\[
\tilde\bfdelta{}' \tilde\bfSigma{}^{-1} \bfSigma\tilde\bfSigma{}^{-1}
\tilde\bfdelta= \tilde\bfdelta{}' \tilde\bfSigma{}^{-1}
\tilde\bfdelta[1+O_P(d_n)]
= \tilde\bfdelta{}' \bfSigma^{-1} \tilde\bfdelta[1+O_P(d_n)] .
\]
Without loss of generality, we assume that
$\tilde\bfdelta= ( \tilde\bfdelta{}_1' , \bfzero' )'$, where
$\tilde{\bfdelta}_1$ is the $\hat{q}$-vector containing nonzero components
of $\tilde\bfdelta$. Let $\bfdelta= ( \bfdelta_1' , \bfdelta_0'
)'$, where
$\bfdelta_1$ has dimension $\hat{q}$.
From Lemma \ref{lem2}(ii),
$\| \tilde\bfdelta_1 - \bfdelta_1 \|^2 =O_P(q_n/n)$ and,
with probability tending to 1,
\[
\| \bfdelta_0\|^2
= \sum_{j: | \hat\delta_j| \leq a_n} \delta_j^2
\leq\sum_{j: | \delta_j| \leq r a_n} \delta_j^2
\leq(r a_n)^{2(1-g)} \sum_{j: | \delta_j| \leq r a_n} \delta_j^{2g}
= O\bigl( a_n^{2(1-g)} D_{g,p}\bigr).
\]
Let $k_n = \max\{ a_n^{2(1-g)} D_{g,p} , q_n/n\}$. Then
$ \| \tilde\bfdelta- \bfdelta\|^2 =
\| \tilde\bfdelta_1 - \bfdelta_1 \|^2 + \|\bfdelta_0\|^2 = O_P(k_n) $.
This together with (\ref{conds})--(\ref{condd})
implies that $(\tilde\bfdelta- \bfdelta)' \bfSigma^{-1}
(\tilde\bfdelta- \bfdelta)= O_P(k_n), $ and hence
\begin{eqnarray*}
\tilde\bfdelta{}' \bfSigma^{-1} \tilde\bfdelta
& = &
\zetap+ 2 \bfdelta' \bfSigma( \tilde\bfdelta- \bfdelta)+
( \tilde\bfdelta- \bfdelta)' \bfSigma^{-1}
( \tilde\bfdelta- \bfdelta)\\
& = & \zetap\bigl[ 1+ O_P \bigl(\sqrt{k_n}/\zetaps\bigr) +
O_P(k_n/\zetap) \bigr]\\
& = & \zetap\bigl[ 1+ O_P \bigl(\sqrt{k_n}/\zetaps\bigr)
\bigr].
\end{eqnarray*}
Write
\begin{eqnarray*}
\bfSigma&=& \pmatrix{
\bfSigma_1 & \bfSigma_{12} \cr
\bfSigma_{12}' & \bfSigma_2 } ,
\qquad
\bfSigma^{-1} = \pmatrix{
\bfC_1 & \bfC_{12} \cr
\bfC_{12}' & \bfC_2 } ,
\\
\tilde\bfSigma &=& \pmatrix{
\tilde\bfSigma_1 & \tilde\bfSigma_{12}\vspace*{1pt} \cr
\tilde\bfSigma_{12}' & \tilde\bfSigma_2 } ,
\qquad
\tilde\bfSigma{}^{-1} = \pmatrix{
\tilde\bfC_1 & \tilde\bfC_{12}\vspace*{1pt} \cr
\tilde\bfC_{12}' & \tilde\bfC_2 } ,
\end{eqnarray*}
where $\bfSigma_1$, $\tilde\bfSigma_1$, $\bfC_1$ and $\tilde\bfC
_1$ are
$q_n \times q_n$ matrices with $q_n$ defined in Lemma~\ref{lem2}(ii). Then
\[
\bfC_{12} = - \bfSigma_1^{-1} \bfSigma_{12} \bfC_2 \quad \mbox
{and}\quad
\tilde\bfC_{12} = - \tilde\bfSigma_1^{-1} \tilde\bfSigma_{12}
\tilde\bfC_2 .
\]
If $\check\bfdelta_1 = ( \tilde\bfdelta{}_1' , \bfzero' )'$ and
$\bar{\mathbf{x}}_1 - \bfmu_1 = ( \bfxi_1' , \bfxi_0' )'$,
where $\check\bfdelta_1$ and $\bfxi_1$ have dimension $q_n$, then
\[
\tilde\bfdelta{}' \tilde\bfSigma{}^{-1} ( \bar{\mathbf{x}}_1 - \bfmu
_1 )
= \check\bfdelta{}_1' \tilde{\bfC}_1 \bfxi_1 + \check\bfdelta{}_1'
\tilde\bfC_{12}
\bfxi_0 = \check\bfdelta{}_1' \tilde{\bfC}_1 \bfxi_1 - \check
\bfdelta{}_1'
\tilde\bfSigma{}_1^{-1} \tilde\bfSigma_{12} \tilde\bfC_2 \bfxi_0 .
\]
Since $\bfxi_1$ has dimension $q_n$,
\[
( \check\bfdelta{}_1' \tilde\bfC_1 \bfxi_1 )^2
\leq(\bfxi_1' \tilde\bfC_1 \bfxi_1 ) (\check\bfdelta{}_1' \tilde
\bfC_1 \check\bfdelta_1)
= (\bfxi_1' \tilde\bfC_1 \bfxi_1 )(\tilde\bfdelta{}' \tilde
\bfSigma{}^{-1} \tilde\bfdelta)
= O_P(q_n /n )(\tilde\bfdelta{}' \tilde\bfSigma{}^{-1} \tilde\bfdelta)
\]
and hence
\[
\check\bfdelta{}_1' \tilde\bfC_1 \bfxi_1
= O_P \bigl(\sqrt{k_n} \bigr) \sqrt{\tilde\bfdelta{}' \tilde
\bfSigma{}^{-1} \tilde\bfdelta}.
\]
Since $\tilde\bfSigma{}_1^{-1} \leq\tilde\bfC_1$,
\begin{eqnarray*}
(\check\bfdelta{}_1' \tilde\bfSigma{}_1^{-1} \tilde\bfSigma_{12}
\tilde\bfC_2 \bfxi_0)^2
& \leq& ( \check\bfdelta{}_1'\tilde\bfSigma{}_1^{-1}
\check\bfdelta_1 )
( \bfxi_0' \tilde\bfC_2 \tilde\bfSigma{}_{12}'
\tilde\bfSigma{}_1^{-1} \tilde\bfSigma_{12} \tilde\bfC_2 \bfxi_0 )
\\
&\leq& ( \tilde\bfdelta{}_1'\tilde\bfC_1 \tilde\bfdelta_1 )
( \bfxi_0' \tilde\bfC_2 \tilde\bfSigma{}_{12}'
\tilde\bfSigma{}_1^{-1} \tilde\bfSigma_{12} \tilde\bfC_2 \bfxi_0 )
\\
&= &
( \tilde\bfdelta{}' \tilde\bfSigma{}^{-1} \tilde\bfdelta)
( \bfxi_0' \tilde\bfC_2 \tilde\bfSigma{}_{12}'
\tilde\bfSigma{}_1^{-1} \tilde\bfSigma_{12} \tilde\bfC_2 \bfxi_0 ).
\end{eqnarray*}
From result (\ref{bl}),
\[
\bfxi_0' \tilde\bfC_2 \tilde\bfSigma{}_{12}'
\tilde\bfSigma{}_1^{-1} \tilde\bfSigma_{12} \tilde\bfC_2 \bfxi_0
= \bfxi_0' \bfC_2 \bfSigma_{12}'\bfSigma_1^{-1} \bfSigma_{12} \bfC
_2 \bfxi_0
[1+O_P(d_n)] .
\]
Under condition (\ref{conds}), all eigenvalues of sub-matrices of
$\bfSigma$ and $\bfSigma^{-1}$ are bounded by $c_0$.
Repeatedly using condition (\ref{conds}), we obtain that
\begin{eqnarray*}
E ( \bfxi_0' \bfC_2 \bfSigma_{12}'\bfSigma_1^{-1} \bfSigma_{12}
\bfC_2 \bfxi_0 )
& \leq& c_0 E ( \bfxi_0' \bfC_2 \bfSigma_{12}' \bfSigma_{12}
\bfC_2 \bfxi_0 ) \\
& = & c_0 n^{-1} \operatorname{trace} (
\bfSigma_{12} \bfC_2 \bfSigma_2 \bfC_2 \bfSigma_{12}' ) \\
& \leq& c_0^4 n^{-1} \operatorname{trace} (
\bfSigma_{12} \bfSigma_{12}' ) \\
& = &
\frac{c_0^4}{n} \sum_{j=1}^{q_n} \sum_{l=q_n+1}^p \sigma_{jl}^2 \\
& \leq& \frac{c_0^{6-h} q_n}{n}
\max_{l \leq p} \sum_{j=1}^p |\sigma_{jl}|^h \\
& = & O( C_{h,p} q_n /n ),
\end{eqnarray*}
where $h$ and $C_{h,p}$ are given in (\ref{cp}). This proves that
\[
\frac{\tilde\bfdelta{}' \tilde\bfSigma^{-1} ( \bar\bfx_1 - \bfmu_1 ) }{\sqrt{
\tilde\bfdelta {}' \tilde\bfSigma^{-1}\bfSigma \tilde\bfSigma^{-1} \tilde\bfdelta}}
= \frac{O_P(\sqrt{k_n}) + O_P( \sqrt{C_{h,p} q_n/n} )}{
\sqrt{1+O_P(d_n)} },
\]
which also holds when $ \bar\bfx_1 - \bfmu_1 $ is replaced by
$ \bar\bfx_2 - \bfmu_2 $ or $\hat\bfdelta - \bfdelta$. Note that
\begin{eqnarray*}
 \hat{\bfdelta}{}' \tilde{\bfSigma}^{-1} \tilde{\bfdelta}
& = & \tilde{\bfdelta}{}' \tilde{\bfSigma}^{-1} \tilde{\bfdelta}
+  ( \hat{\bfdelta} - \bfdelta ) ' \tilde{\bfSigma}^{-1} \tilde{\bfdelta}
+  ( \bfdelta - \tilde{\bfdelta} ) ' \tilde{\bfSigma}^{-1} \tilde{\bfdelta} \\
& = & \tilde{\bfdelta}{}' \tilde{\bfSigma}^{-1} \tilde{\bfdelta}
+  ( \hat{\bfdelta} - \bfdelta ) ' \tilde{\bfSigma}^{-1} \tilde{\bfdelta}
+  \Delta_p O_P \bigl( \sqrt{k_n} \bigr).
\end{eqnarray*}
Therefore,
\begin{eqnarray*}
\frac{(-1)^k \tilde\bfdelta{}' \tilde\bfSigma^{-1} ( \bfmu_k - \bar\bfx_k ) -
\hat\bfdelta {}' \tilde\bfSigma^{-1} \tilde\bfdelta /2}{\sqrt{
\tilde\bfdelta {}' \tilde\bfSigma^{-1}\bfSigma \tilde\bfSigma^{-1}
\tilde\bfdelta}}&=&
\frac{O_P(\sqrt{k_n}) + O_P( \sqrt{C_{h,p} q_n/n} )}{
\sqrt{1+O_P(d_n)} } \\
& &{} - \frac{\zetaps\sqrt{1+O_P(\sqrt{k_n}/\zetaps)}}{
2 \sqrt{1+O_P(d_n)}} \\
& = & O_P \bigl(\sqrt{k_n} \bigr) + O_P \bigl( \sqrt
{C_{h,p}q_n /n} \bigr) \\
& &{} - \frac{\zetaps}{2} \bigl[ 1+ O_P \bigl(\sqrt
{k_n}/\zetaps\bigr) +O_P(d_n) \bigr] \\
& = & - \frac{\zetaps}{2} \biggl[ 1 + O_P \biggl( \frac{\sqrt{
C_{h,p}q_n }}{\zetaps\sqrt{n}} \biggr) \\
& &\hphantom{- \frac{\zetaps}{2} \biggl[}{} + O_P \biggl(\frac{\sqrt{k_n}}{\zetaps} \biggr)
+O_P(d_n) \biggr]\\
& = & - \frac{\zetaps}{2} [ 1+O_P (b_n
) ].
\end{eqnarray*}
%
%
This proves the result in (i).
The proofs of (ii)--(iv) are the same as the proofs
for Theorem \ref{thm1}(ii)--(iv) with $s_n$ replaced by $b_n$.
This completes the proof.
\end{pf*}

\section*{Acknowledgments}

The authors would like to
thank two referees and an associate editor
for their helpful comments and suggestions, and Dr. Weidong Liu
for his help in correcting an error in the proof of Theorem 3.


%

\printaddresses

\end{document}